\documentclass[conference]{IEEEtran}
\usepackage{cite}

\ifCLASSINFOpdf
  \usepackage[pdftex]{graphicx}
  \graphicspath{figs/}
  \DeclareGraphicsExtensions{.pdf,.eps,.png}
\else
\fi
\usepackage[cmex10]{amsmath}
\usepackage{amssymb}
\usepackage{algorithmic}

\usepackage{array}

\usepackage{mdwmath}
\usepackage{mdwtab}

\usepackage{eqparbox}

\usepackage[font=footnotesize,caption=false]{subfig}
\usepackage{url}

\graphicspath{{./}{./figs/}}

  \def\clap#1{\hbox to 0pt{\hss#1\hss}}

\usepackage{bm}
\providecommand{\mat}[1]{\bm{#1}}%
\renewcommand{\vec}[1]{\mathbf{#1}}

\providecommand{\mA}{\ensuremath{\mat{A}}}
\providecommand{\mB}{\ensuremath{\mat{B}}}
\providecommand{\mC}{\ensuremath{\mat{C}}}

\providecommand{\mE}{\ensuremath{\mat{E}}}

\providecommand{\mG}{\ensuremath{\mat{G}}}
\providecommand{\mH}{\ensuremath{\mat{H}}}

\providecommand{\mW}{\ensuremath{\mat{W}}}

\providecommand{\va}{\ensuremath{\vec{a}}}

\providecommand{\vm}{\ensuremath{\vec{m}}}

\providecommand{\vs}{\ensuremath{\vec{s}}}

\providecommand{\vv}{\ensuremath{\vec{v}}}
\providecommand{\vw}{\ensuremath{\vec{w}}}
\providecommand{\vx}{\ensuremath{\vec{x}}}

\providecommand{\vz}{\ensuremath{\vec{z}}}

\newcommand{\hmC}{\widehat{\mC}}
\newcommand{\hmG}{\widehat{\mG}}
\newcommand{\hLambda}{\widehat{\Lambda}}

\newcommand{\hmW}{\widehat{\mW}}
\newcommand{\hlambda}{\hat{\lambda}}

\newcommand{\sM}{\mathcal{M}}

\newcommand{\sP}{\mathcal{P}}

\newcommand{\bmat}[1]{\begin{bmatrix}#1\end{bmatrix}}

\begin{document}
%
\title{Computing Active Subspaces Efficiently with Gradient Sketching}

\author{\IEEEauthorblockN{Paul G.~Constantine}
\IEEEauthorblockA{Applied Mathematics and Statistics\\
Colorado School of Mines\\
Golden, Colorado 80401\\
Email: \url{pconstan@mines.edu}}
\and
\IEEEauthorblockN{Armin Eftekhari and Michael B.~Wakin}
\IEEEauthorblockA{Electrical Engineering and Computer Science\\
Colorado School of Mines\\
Golden, Colorado 80401\\
Email: \url{aeftekha,mwakin@mines.edu}}
}


%


\maketitle

\begin{abstract}
Active subspaces are an emerging set of tools for identifying and exploiting the most important directions in the space of a computer simulation's input parameters; these directions depend on the simulation's quantity of interest, which we treat as a function from inputs to outputs. To identify a function's active subspace, one must compute the eigenpairs of a matrix derived from the function's gradient, which presents challenges when the gradient is not available as a subroutine. We numerically study two methods for estimating the necessary eigenpairs using only linear measurements of the function's gradient. In practice, these measurements can be estimated by finite differences using only two function evaluations, regardless of the dimension of the function's input space. 
\end{abstract}


%
\IEEEpeerreviewmaketitle

\section{Active subspaces}

\noindent Modern physics and engineering simulations take several inputs---e.g., boundary conditions, material properties, and forcings---and output several quantities of interests. The scientist uses these simulations to study the relationship between inputs and outputs. Uncertainty quantification seeks precise characterization of the simulation's quantities of interest subject to variability in the inputs. These characterizations often reduce to parameter studies---such as optimization, integration, or response surface modeling---that treat the simulation as a mapping between inputs $\vx$ and a quantity of interest $f(\vx)$. However, thorough parameter studies quickly become infeasible as the dimension of $\vx$ grows, particularly if evaluating $f(\vx)$ (i.e., running the physical simulation) is computationally expensive. To combat this curse of dimensionality, one may seek a low-dimensional parameterization of $f(\vx)$ that (i) maintains the input/output representation and (ii) enables otherwise infeasible parameter studies. One idea is to identify the least important input parameters and fix them at nominal values, thus reducing the dimension of the parameter study. Such identification is the domain of \emph{sensitivity analysis}, and several techniques exist that use a few simulation runs to screen the inputs' importance---such as local perturbations, elementary effects, or sensitivity indices~\cite{saltelli2008global}. A more general approach is to identify important linear combinations of the inputs $\vx$ and focus parameter studies along the associated directions. \emph{Active subspaces} are defined by important directions in the high-dimensional space of inputs; once identified, the scientist can exploit the active subspace to enable otherwise infeasible parameter studies for expensive simulations~\cite{asm2015}. 

Assume that $\vx\in\mathbb{R}^m$ is a vector of simulation inputs, and let the input space be equipped with a probability density function $\rho(\vx)$ that is strictly positive in the domain of $f$ and zero outside the domain. In practice, $\rho$ identifies the set of inputs of interest and quantifies the variability. We assume that the independent inputs have been shifted and scaled to be centered at the origin and have equal variances. Assume that $f:\mathbb{R}^m\rightarrow\mathbb{R}$ is continuous, square-integrable with respect to $\rho$, and differentiable with gradient vector $\nabla f\in\mathbb{R}^m$; also assume that $f$'s gradient is square-integrable with respect to $\rho$. The active subspace is defined by the first $n<m$ eigenvectors of the following $m\times m$ symmetric positive semi-definite matrix,
\begin{equation}
\label{eq:cmat}
\mC \;=\; \int \nabla f\,\nabla f^T\,\rho\,d\vx \;=\; \mW\Lambda\mW^T, 
\end{equation}
where the non-negative eigenvalues are ordered in descending order. The eigenvalue $\lambda_i$ measures the average change in $f$ subject to perturbations in $\vx$ along the corresponding eigenvector $\vw_i$, 
\begin{equation}
\lambda_i \;=\; \int (\nabla f^T\vw_i)^2\,\rho\,d\vx.
\end{equation}
For example, if $\lambda_i=0$, then $f$ is constant along the direction $\vw_i$. If $f$ is constant along a direction, then one can ignore this direction when studying the behavior of $f$ under changes in $\vx$. Suppose that the first $n<m$ eigenvalues are much larger than the trailing $m-n$, and let $\mW_1$ be the first $n$ columns of the orthogonal eigenvector matrix $\mW$. Then a reasonable approximation for $f$ is
\begin{equation}
\label{eq:gmodel}
f(\vx) \;\approx\; g(\mW_1^T\vx),
\end{equation} 
where $g$ is a properly constructed map from $\mathbb{R}^n$ to $\mathbb{R}$. In~\cite{Constantine2014}, we study a particular construction of $g$, its approximation properties, and how those properties change when $\mW_1$ is estimated. 

There are several recent works that study models similar to \eqref{eq:gmodel} in uncertainty quantification~\cite{Tipireddy2014,Stoyanov2015}, computational engineering~\cite{Bang2012,Berguin2014a,constantine2014exploiting,Russi2010}, and approximation theory~\cite{Fornasier2012,Tyagi2014}. In statistics, subspace-based dimension reduction has wide use in regression modeling under the monikers \emph{sufficient dimension reduction}~\cite{cook2009regression} and \emph{efficient dimension reduction}~\cite{Li1992}. Many in that community have recognized the connection between the dimension reduction space and the matrix $\mC$ in  \eqref{eq:cmat}; see~\cite{samarov1993exploring,hristache2001}. However, in contrast to our case where the gradient $\nabla f$ is given, the gradient of the regression function with respect to predictors must be learned. Additionally, there is no random noise in computer simulations in contrast to regression modeling. 

The tremendous potential benefits of dimension reduction drive us to pursue methods to estimate the eigenvalues $\Lambda$ and eigenvectors $\mW$ from \eqref{eq:cmat}. In~\cite{constantine2014computing}, we analyze the following Monte Carlo method. First, draw $M$ samples $\{\vx_i\}$ independently according to $\rho$. For each $\vx_i$, compute $\nabla f_i = \nabla f(\vx_i)$. Then compute
\begin{equation}
\label{eq:capprox}
\mC \;\approx\; \hmC \;=\; \frac{1}{M}\sum_{i=1}^M \nabla f_i\,\nabla f_i^T
\;=\; \hmW\hLambda\hmW^T.
\end{equation}
In~\cite{constantine2014computing}, we use tools from nonasymptotic random matrix theory to study the approximation error in the estimated eigenvalues and subspaces~\cite{gittens2011tail,tropp2012user}. The approach assumes that one has access to the gradient $\nabla f(\vx)$ as a black box, which is not unreasonable in modern simulation codes due to adjoint methods or algorithmic differentiation~\cite{Griewank00}. However, many legacy simulation codes lack gradient routines, and one must estimate the eigenpairs with only evaluations of $f(\vx)$.

\section{Estimating active subspaces without gradients}
\noindent When the gradient is not available as a black box, one might estimate the gradient through evaluations of $f$, e.g., with finite differences or related techniques~\cite{conn2009introduction}. A first-order finite difference approximation to the gradient vector requires $m+1$ evaluations of $f$, so estimating $\hmC$ in \eqref{eq:capprox} takes $M(m+1)$ evaluations, which may be prohibitively expensive for large-scale simulations. To surmount this challenge, we take advantage of a finite difference approximation of directional derivatives. For sufficiently small $h>0$ and $\va\in\mathbb{R}^m$,
\begin{equation}
\label{eq:fddir}
\nabla f(\vx)^T\va \;\approx\; (f(\vx+h\va)-f(\vx))/h.
\end{equation}
The left side of \eqref{eq:fddir} can be interpreted as a linear measurement of the gradient vector, and the right side takes only two evaluations of $f$---regardless of the dimension $m$. We study two techniques that exploit \eqref{eq:fddir} to estimate $\mW_1$ from \eqref{eq:cmat} using only evaluations of $f$. The relationship \eqref{eq:fddir} is also exploited in the works~\cite{Fornasier2012,Tyagi2014} to reconstruct functions of the form $f(\vx)=g(\mA\vx)$, where $\mA\in\mathbb{R}^{n\times m}$. Note the subtle difference between this reconstruction and the approximation model in \eqref{eq:gmodel}. Our goal in this paper is to estimate $\mW_1$, not propose a choice for $g$ in \eqref{eq:gmodel}. Another important difference is in the linear measurement operator; \cite{Fornasier2012,Tyagi2014} use the same operator for each $\nabla f_i$ in \eqref{eq:capprox}, while we use an independent measurement operator for each $\nabla f_i$.

\subsection{Eigenvectors from projections}
\label{sec:qi}
\noindent Our first approach is based on work by Qi and Hughes~\cite{Qi2012} for estimating principal components from linear measurements of a collection of vectors. Suppose that $\vz_i\in\mathbb{R}^m$, $i=1,\dots,M$, are independent, zero-mean random vectors of the form
\begin{equation}
\label{eq:zmodel}
\vz_i \;=\; \sum_{j=1}^d w_{ij}\,\sigma_j\,\vv_j, \qquad w_{ij}\sim\mathcal{N}(0,1),
\end{equation}
where the $\vv_j$ are orthonormal vectors. Let $\mE_i\in\mathbb{R}^{m\times k}$, $k<m$, have independent standard Gaussian entries, and define the measurements $\vm_i \;=\; \mE_i^T\vz_i$. Then the orthogonal projection of $\vz_i$ onto the column space of $\mE_i$ is
\begin{equation}
\label{eq:projz}
\sP_i\vz_i \;:=\; \mE_i(\mE_i^T\mE_i)^{-1}\vm_i \;=\;
\mE_i(\mE_i^T\mE_i)^{-1}\mE_i^T\vz_i.
\end{equation}
Qi and Hughes~\cite[Theorem 2]{Qi2012} show that the first $d$ eigenvectors of the matrix
\begin{equation}
\frac{1}{M} \sum_{i=1}^M \sP_i\vz_i\,(\sP_i\vz_i)^T
\end{equation}
converge to the vectors $\vv_1,\dots,\vv_d$ as $M$ goes to infinity. The eigenvalues converge to quantities related to $\sigma_i$ from \eqref{eq:zmodel}. In fact, their results show how to estimate the mean and principal components even when the model \eqref{eq:zmodel} has a nonzero mean and random noise.

There is no reason to suspect that the gradient vectors $\nabla f_i$ from \eqref{eq:capprox} satisfy the model \eqref{eq:zmodel}. Nonetheless, we can numerically check if the eigenvectors of the matrix
\begin{equation}
\label{eq:cproj}
\hmC_{\sP} \;=\; \frac{1}{M} \sum_{i=1}^M \sP_i\nabla f_i\,(\sP_i\nabla f_i)^T,
\end{equation}
where $\sP_i$ is defined as in \eqref{eq:projz}, are close to the eigenvectors of \eqref{eq:capprox} for chosen test problems. If so, we can exploit the finite difference relationship \eqref{eq:fddir} to efficiently estimate each element of the $k$-vector $\mE_i^T\nabla f_i$, which is the analogue of $\vm_i$ in \eqref{eq:projz}. In this case, we can estimate the eigenvectors $\hmW$ from \eqref{eq:capprox} using $M(k+1)<M(m+1)$ evaluations of $f$. Note that, by the analysis in~\cite{Qi2012}, we do not expect the eigenvalues of $\hmC_{\sP}$ to converge to those of $\hmC$. 

\subsection{Low-rank approximation from linear measurements}
\label{sec:altmin}
\noindent Our second approach is based on low-rank approximation of the matrix of gradients $\hmG = [\nabla f_1,\dots,\nabla f_M]\in\mathbb{R}^{m\times M}$ using linear measurements of the gradients. Define the linear measurement operator $\sM(\cdot)$ as
\begin{equation}
\sM(\hmG) \;:=\; \bmat{\mE_1^T\nabla f_1 & \cdots & \mE_M^T\nabla f_M},
\end{equation}
where $\mE_i\in\mathbb{R}^{m\times k}$ has independent standard Gaussian entries as in the previous section. Let $r$ be the rank of the low-rank approximation. We seek matrices $\mA\in\mathbb{R}^{m\times r}$ and $\mB\in\mathbb{R}^{M\times r}$ that solve
\begin{equation}
\label{eq:altmin}
\underset{\mA,\,\mB}{\operatorname{minimize}}\; \|\sM(\hmG) - \sM(\mA\mB^T)\|_F,
\end{equation}
where $\|\cdot\|_F$ is the Frobenius norm. We estimate the minimizers with alternating least-squares. Given $\mA$, \eqref{eq:altmin} is a linear least-squares problem for $\mB$. Similarly, given $\mB$, \eqref{eq:altmin} is a linear least-squares problem for $\mA$. We choose a starting value for $\mA$ as the first $r$ eigenvectors, scaled by the square-roots of the first $r$ eigenvalues, of $\hmC_{\sP}$ in \eqref{eq:cproj}. In this sense, the alternating least-squares can be considered an iterative refinement on the estimates from the first method. Once $\mA$ and $\mB$ are estimated, we compute the left singular vectors of $\mA\mB^T$ to estimate $\mW_1$ in \eqref{eq:gmodel}. 

The low-rank model $\mA\mB^T$ requires the user to choose the rank $r$. To solve the least-squares subproblem for $\mA$ without additional regularization, we need $r$ less than the number $k$ of linear measurements (i.e., the number of columns in $\mE_i$). Additionally, we are guided by the ultimate goal of the dimension reduction, which is to construct $g$ in \eqref{eq:gmodel}. Without prior knowledge of the relationship between $f$ and $\vx$, constructing a response surface is generally too expensive in more than a handful of dimensions. Thus, we may reasonably keep $r$ less than 8 or 9. If there is no gap within the first 8 or 9 eigenvalue estimates, then \eqref{eq:gmodel} may be inappropriate.

\section{Experiments}
\noindent We test the two methods on two functions: (i) a quadratic polynomial in $m=10$ dimensions and (ii) a quantity of interest from the solution of a PDE whose operator coefficients depend on $m=100$ variables. In each experiment, we compute $\hmC$ from \eqref{eq:capprox} with a fixed number $M$ of gradient samples, and we consider the eigenpairs $\hmW$ and $\hLambda$ to be the \emph{true} values. We compute the error in the first six eigenvalue estimates and the error in the subspaces defined by the eigenvectors estimates as
\begin{equation}
\label{eq:evalerr}
\left(\frac{\sum_{i=1}^6 (\hlambda_i-\tilde{\lambda}_i)^2}{\sum_{i=1}^6 \hlambda_i^2}\right)^{1/2}
\end{equation}
and
\begin{equation}
\label{eq:suberr}
\|\hmW_1\hmW_1^T - \tilde{\mW}_1\tilde{\mW}_1^T\|_2,
\end{equation}
respectively, where $\tilde{\lambda}_i$ and $\tilde{\mW}$ are the eigenvalue and eigenvector estimates, respectively, from the linear measurement-based approaches. The subscript 1 on the matrices $\hmW$ and $\tilde{\mW}$ indicates that they contain only the first $n$ columns; we note $n$ when needed. We do not expect the eigenvalue estimates from the projection-based method in Section \ref{sec:qi} to converge, but we report them to compare with the errors from the alternating least-squares approach. We repeat the study 20 times with independently drawn Gaussian measurement matrices $\mE_i$. The resulting errors are averaged over these 20 trials. We do not study the approximation properties of the finite differences in \eqref{eq:fddir}, e.g., by varying the finite difference parameter $h$. Instead, we examine cases where the gradient is available, thus focusing on the performance of the linear measurement-based methods when the answer is known. 
 
\subsection{Quadratic function}
\begin{figure*}[!t]
\centerline{
\subfloat[]{\includegraphics[width=0.24\textwidth]{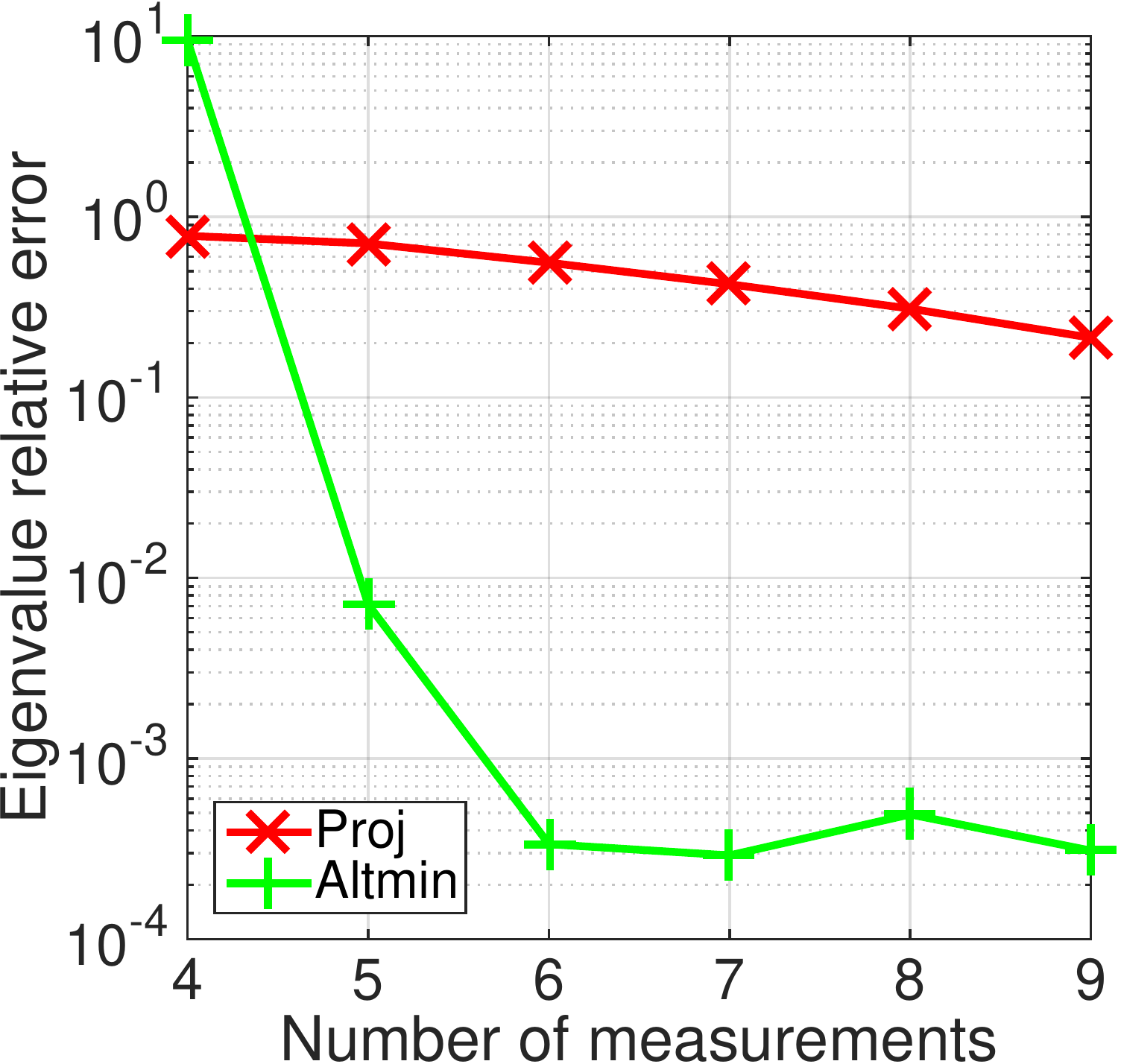}%
\label{fig:eval_err_quad}}
\hfil
\subfloat[]{\includegraphics[width=0.24\textwidth]{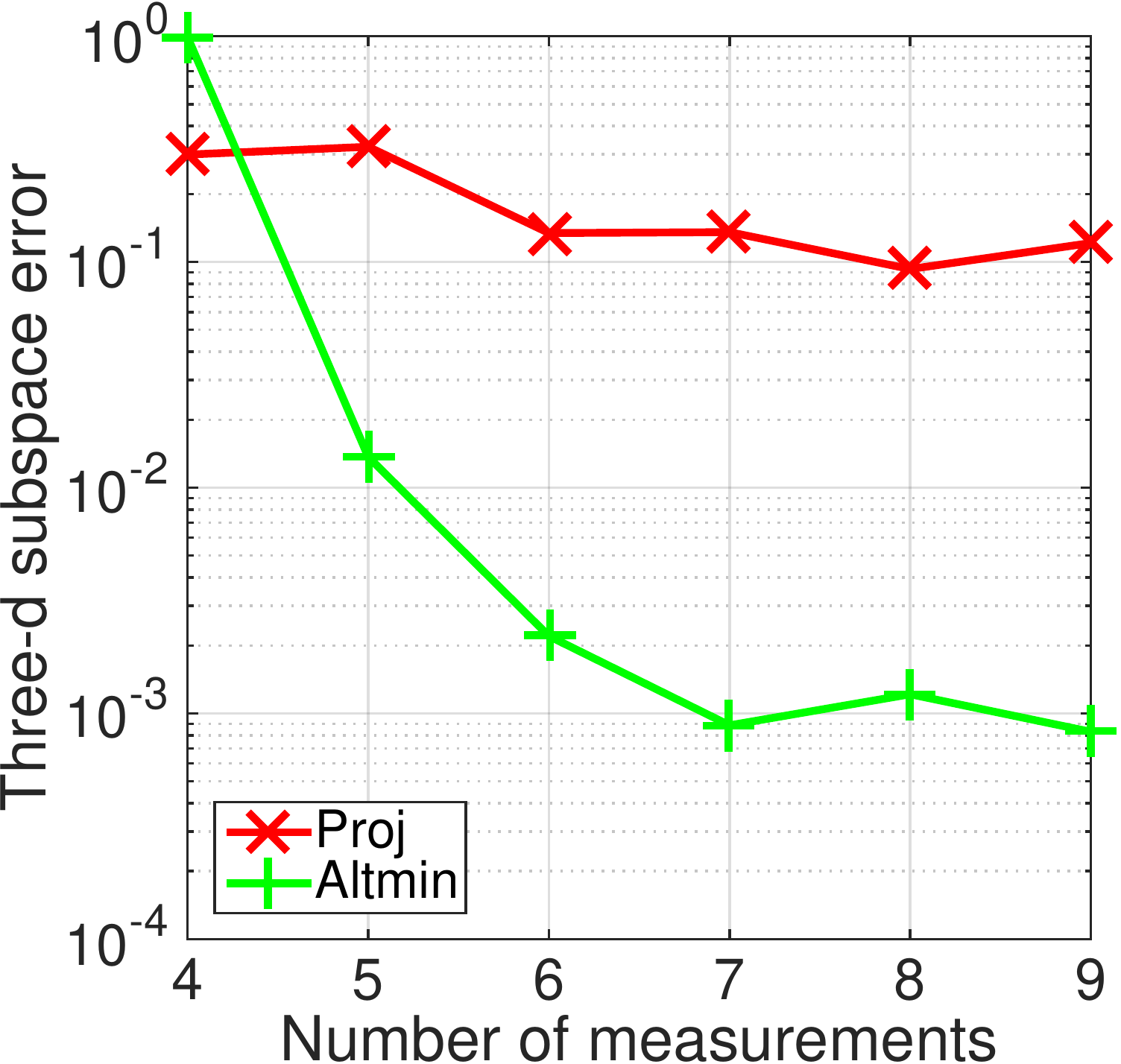}%
\label{fig:sub1_err_quad}}
\hfil
\subfloat[]{\includegraphics[width=0.24\textwidth]{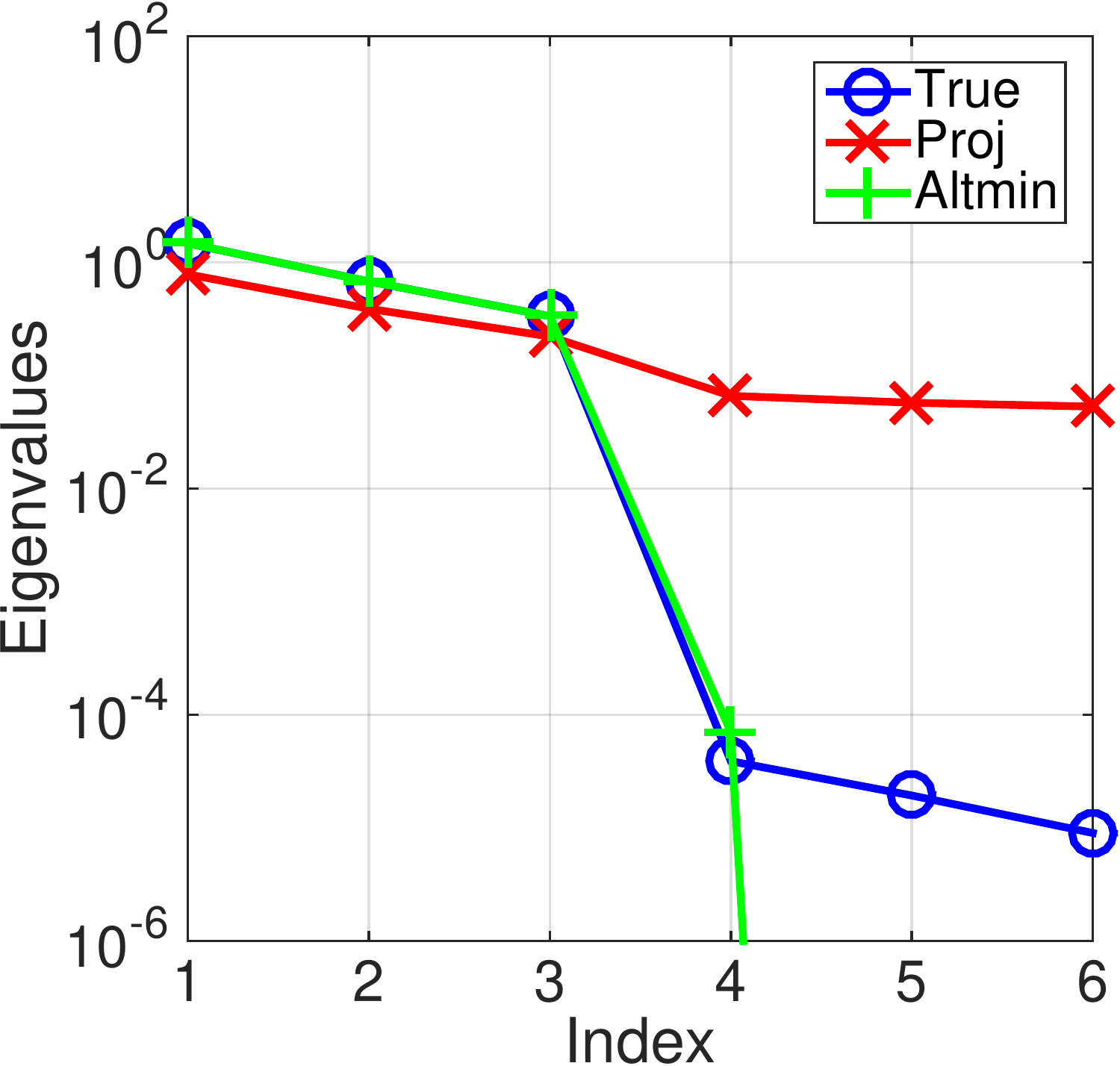}%
\label{fig:evals_quad}}
\hfil
\subfloat[]{\includegraphics[width=0.24\textwidth]{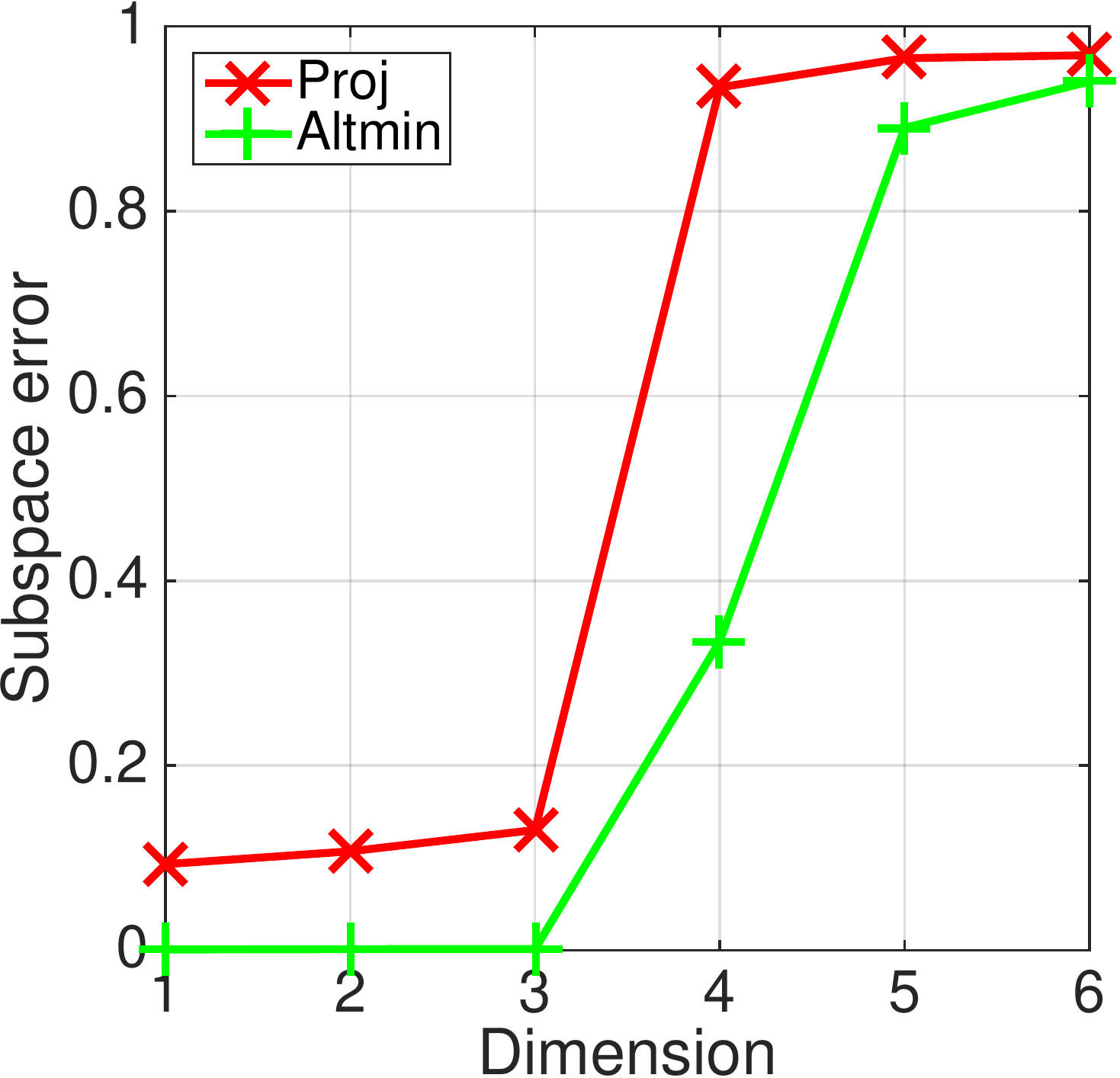}%
\label{fig:sub_quad}}
}
\caption{In the legends, ``Proj'' indicates the projection-based method from Qi and Hughes~\cite{Qi2012} described in Section \ref{sec:qi}, and ``Altmin'' indicates the alternating least-squares minimization from Section \ref{sec:altmin}. Figure \ref{fig:eval_err_quad} shows the error \eqref{eq:evalerr} in the first six eigenvalue estimates as a function of the number $k$ of measurements; note that we do not expect the eigenvalues from the projection-based method to converge. Figure \ref{fig:sub1_err_quad} shows the error in the estimate of the three-dimensional active subspace \eqref{eq:suberr} as a function of the number of measurements. Figure \ref{fig:evals_quad} shows the first six eigenvalues using $k=7$ measurements. Figure \ref{fig:sub_quad} shows the errors in the subspace estimates for dimension $n=1$ to $n=6$ using $k=7$ measurements.}
\label{fig:quad}
\end{figure*}

\noindent Let $\mH\in\mathbb{R}^{10\times 10}$ be symmetric and positive semidefinite, and let $f(\vx) = \frac{1}{2}\vx^T\mH\vx$, defined on the domain $\vx\in[-1,1]^{10}$ with a uniform density $\rho$. The gradient is $\nabla f(\vx)=\mH\vx$. The eigenvectors of $\mC$ from \eqref{eq:cmat} are the eigenvectors of $\mH$, and the eigenvalues of $\mC$ are the eigenvalues of $\mH$, squared and divided by 3. We construct $\mH$ so that its eigenvalues decay at a slow spectral rate, except for a large gap between the third and fourth eigenvalues; this indicates that the active subspace is three-dimensional. We compute $\hmC$ from \eqref{eq:capprox} using $M=200$ gradient samples. We study the quality of the estimates as the number $k$ of measurements goes from 4 to 9. Note that $k=10$ measurements would produce a perfect reconstruction, almost surely. For the alternating least-squares method, we choose the rank $r=4$. The caption in Figure \ref{fig:quad} describes the contents of the subfigures showing the results. In general, we see the least-squares method outperform the projection-based method once the number of measurements exceeds the rank $r=4$.


\subsection{PDE model}
\noindent Let $u=u(\vs,\vx)$ solve the Poisson equation in two spatial dimensions,
\begin{equation}
-\nabla_{\vs}\cdot(a\,\nabla_{\vs} u)\;=\;1,\qquad \vs\in[0,1]^2,
\end{equation}
with homogeneous Dirichlet boundary conditions on the left, top, and bottom of the domain and homogeneous Neumann conditions on the right side of the domain. The spatially varying operator coefficients $a=a(\vs,\vx)$ are parameterized by $m=100$ independent parameters $\vx\in\mathbb{R}^{100}$ via a truncated Karhunen-Loeve expansion of a Gaussian random field with a relatively long correlation length. This implies that $\rho$ is a standard Gaussian density in 100 dimensions. The quantity of interest $f(\vx)$ is the average of the solution $u$ over the right side of the spatial domain. The solution is approximated with a well-resolved finite element method, and the gradient $\nabla f(\vx)$ is computed with an adjoint scheme. More information on this problem can be found in~\cite{Constantine2014}. 

\begin{figure*}[!t]
\centerline{
\subfloat[]{\includegraphics[width=0.24\textwidth]{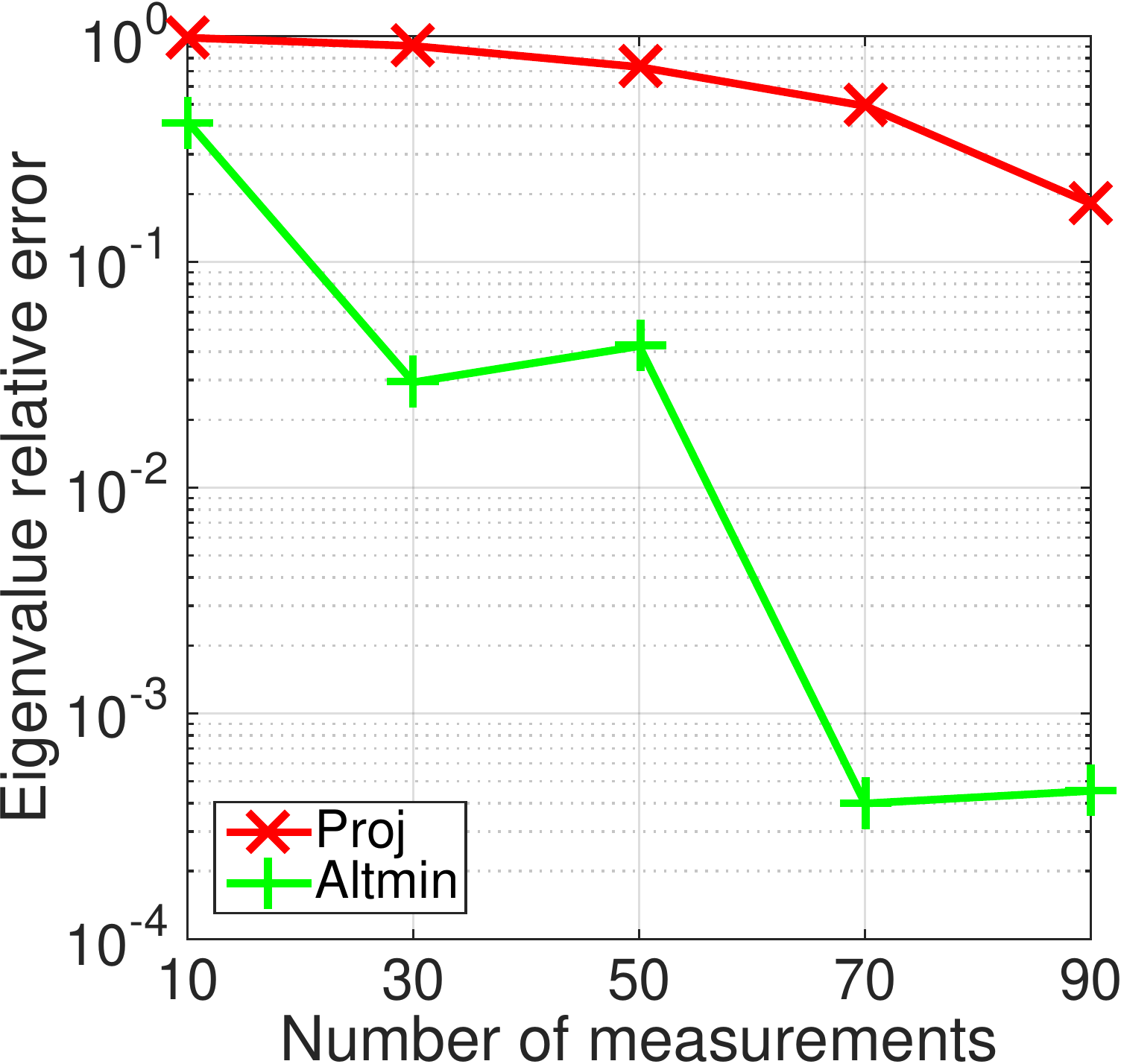}%
\label{fig:eval_err}}
\hfil
\subfloat[]{\includegraphics[width=0.24\textwidth]{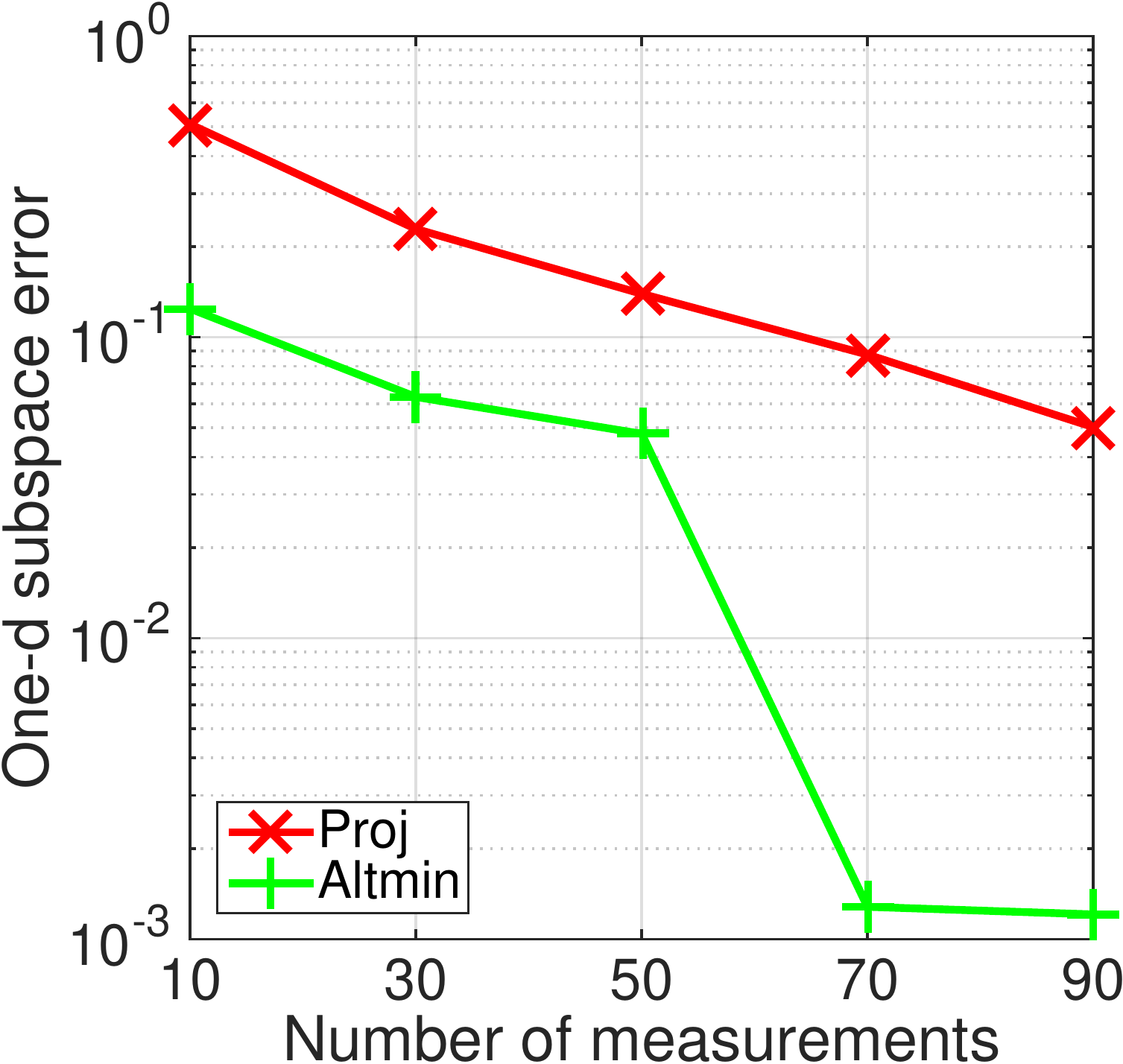}%
\label{fig:sub1_err}}
\hfil
\subfloat[]{\includegraphics[width=0.24\textwidth]{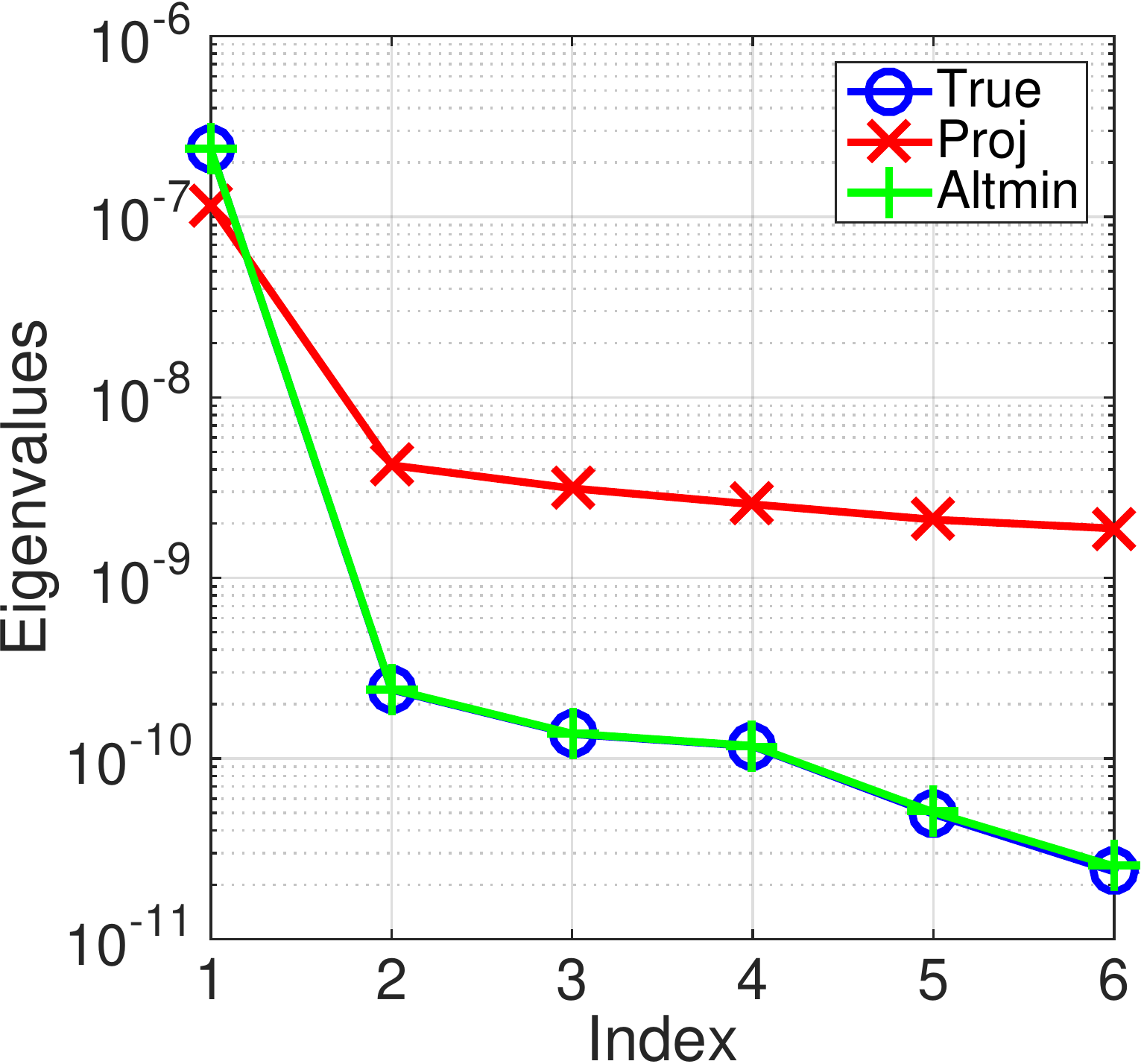}%
\label{fig:evals_pde}}
\hfil
\subfloat[]{\includegraphics[width=0.24\textwidth]{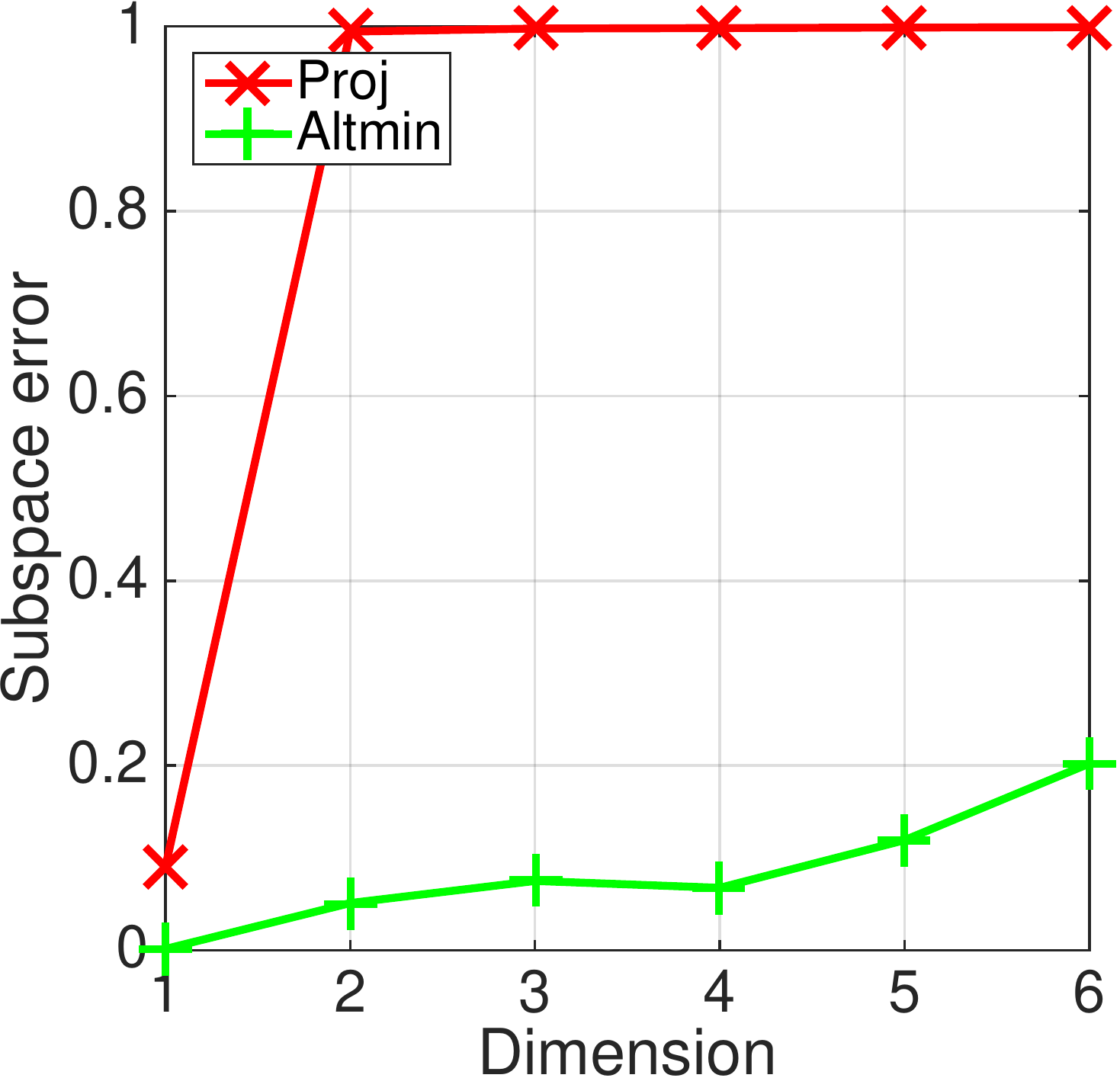}%
\label{fig:sub_pde}}
}
\caption{In the legends, ``Proj'' indicates the projection-based method from Qi and Hughes~\cite{Qi2012} described in Section \ref{sec:qi}, and ``Altmin'' indicates the alternating least-squares minimization from Section \ref{sec:altmin}. Figure \ref{fig:eval_err} shows the error \eqref{eq:evalerr} in the first six eigenvalue estimates as a function of the number $k$ of measurements; note that we do not expect the eigenvalues from the projection-based method to converge. Figure \ref{fig:sub1_err} shows the error in the estimate of the one-dimensional active subspace \eqref{eq:suberr} as a function of the number $k$ of measurements. Figure \ref{fig:evals_pde} shows the first six eigenvalues using $k=70$ measurements. Figure \ref{fig:sub_pde} shows the errors in the subspace estimates for dimension $n=1$ to $n=6$ using $k=70$ measurements.}
\label{fig:pde}
\end{figure*}

The matrix $\hmC$ from \eqref{eq:capprox} is estimated with $M=300$ gradient samples. The particular quantity of interest has a dominant one-dimensional active subspace; in other words, there is a large gap between the first and second eigenvalues of $\hmC$. We study the approximations as the number $k$ of measurements increases by 20 from 10 to 90. For the least-squares method, we choose the rank $r=8$. The caption in Figure \ref{fig:pde} describes the contents of the subfigures showing the results. In general, the alternating least-squares method outperforms the projection-based method. Note the sharp decrease in the error once the number $k$ of measurements exceeds 50. 

\section{Conclusion}
\noindent We have studied the approximation properties of two methods for estimating the eigenvectors and eigenvalues that define a multivariate function's active subspace. These methods use linear measurements of the gradient instead of the full gradient, which can be efficiently computed with finite differences. The first method, based on work by Qi and Hughes~\cite{Qi2012}, uses random projections of the gradient. The second uses a low-rank approximation of the matrix of gradient measurements fit with an alternating least-squares method. The low-rank approximation performs better than the projection-based method in two numerical tests. Future work will (i) analyze the improvement with the low-rank approximation and (ii) study the effect of finite difference approximations of the linear measurements. These initial results are sufficiently promising to pursue such studies. 

\section*{Acknowledgment}
\noindent The first author was partially supported by the U.S. Department of Energy Office of Science, Office of Advanced Scientific Computing Research, Applied Mathematics program under Award Number DE-SC-0011077. The second and third authors were partially supported by NSF CAREER grant CCF-1149225 and NSF grant CCF-1409258.



\bibliographystyle{IEEEtran}
\bibliography{camsap-as}
%
%
%

\end{document}